\documentclass[oneside]{amsart}

\newcommand{\R}{\mathds R}
\newcommand{\I}{\mathds I}

\usepackage{epsfig}
\usepackage{times}
\usepackage{dsfont}
\usepackage{amsmath}
\usepackage[hypertex]{hyperref}

\theoremstyle{plain}\newtheorem{teo}{Theorem}[section]
\theoremstyle{plain}\newtheorem{lem}[teo]{Lemma}
\theoremstyle{plain}\newtheorem{prop}[teo]{Proposition}
\theoremstyle{plain}
\theoremstyle{definition}
\theoremstyle{remark}\newtheorem{rem}[teo]{Remark}
\theoremstyle{remark}\newtheorem{ex}[teo]{Example}
\theoremstyle{plain}

\numberwithin{equation}{section}

\title[Dynamics of homogeneous scalar fields with general
self-interaction potentials]%
{Dynamics of homogeneous scalar fields with general
self-interaction potentials: cosmological and gravitational
collapse models}

\author[R.\ Giamb\`o ,\ F.\ Giannoni]{Roberto Giamb\`o, Fabio Giannoni}
\address{Dipartimento di Matematica e Informatica \hfill\break\indent
Universit\`a di Camerino\hfill\break\indent Italy}
\email{roberto.giambo@unicam.it, fabio.giannoni@unicam.it}

\author[G.\ Magli]{Giulio Magli}
\address{Dipartimento di Matematica \hfill\break\indent Politecnico di
Milano \hfill\break\indent Italy} \email{magli@mate.polimi.it}

\begin{document}
\begin{abstract}

The general relativistic dynamics of a wide class of
self-interacting, self-gravitating homogeneous scalar fields
models is analyzed. The class is characterized by certain general
conditions on the scalar field potential, which include both
asymptotically polynomial and exponential behaviors. Within this
class, we show that the generic evolution is always divergent in a
finite time, and then make use of this result to construct
cosmological models as well as radiating collapsing star models of
the Vaidya type. It turns out that blackholes are generically
formed in such models.

\end{abstract}

\maketitle

\section{Introduction}\label{sec:intro}


Scalar fields as sources of self-gravitating models have attracted
a great deal of attention in cosmology and in relativistic
astrophysics. It is, indeed, worth mentioning the fundamental role
that they play as models of the early universe and, on the other
hand, the relevance of self-gravitating scalar field as test
models for gravitational collapse (see e.g. \cite{Global} and
references therein). There is, however, a crucial difference in
the way scalar fields have been treated in these two scenarios.

First, consider the cosmological models. Here, the scalar field is
always assumed as {\it self-interacting}; in other words, the
cosmological expansion is "driven" by a non-vanishing scalar field
potential, and the "free" case (i.e. the case of the quadratic
field lagrangian) corresponds, as is well known, just to the case
of a dynamical cosmological constant. Thus, in cosmology, the
presence of a non vanishing field potential is a key point, and
indeed many efforts have been made to study the possible dynamical
behaviors of the models in dependence of the choice of the
potential, and even to try to put in evidence possible large-scale
observable effects (see e.g.\cite{st1,st2} and references
therein). The situation is even more intriguing when issues from
string theories come into play; indeed, here one is lead to test
simple models which are asymptotically {\it anti}-De Sitter, and
the potential for the scalar field might have, as a consequence,
"non-standard" behaviors (we shall come back on this interesting
issue later on). All in all, self-interacting fields with a non
trivial potential are a key ingredient in cosmological models and
therefore, of course, in the nature of the cosmological
singularities.

Consider, now, the "astrophysical counterpart" of the cosmological
solutions, that is, gravitational collapse. Here, as is well
known, there still exist, unsolved, the issue of the cosmic
censorship, that is whether the singularities which form at the
end state of the collapse are always covered by an event horizon,
or not \cite{ns,j-rev}. The scalar field model of course proposes itself a good
test-bed also in this different scenario. However, in this context
the scalar fields sources have usually been considered as {\it
free}, i.e. minimally coupled to gravity via equivalence principle
only. For this very special case the existence of naked
singularities has been shown, but it has been also shown that such
singularities are in a sense non generic with respect to the
choice of the initial data \cite{c1,c2}. Only in recent years, a
few results have been added to this scenario with the study of
gravitational collapse of {\it homogeneous} self-interacting
scalar fields. In these papers formation of naked singularities
has been found \cite{hsf,joshi}, but with the choice of very
special potentials.

All in all, in both applications (cosmology, and gravitational
collapse) it would be especially welcome a treatment of
self-interacting scalar field dynamics able to treat the models in
a unified manner, and to predict the qualitative behavior in
dependence of the choice of the potential, viewed as an element of
a space of admissible "equations of state" for the matter source.
It is the aim of the present paper to provide this unified
framework, at least in the special case of {\it homogeneous}
scalar fields. Indeed, we study here the dynamical behavior of
such fields for a very general class of potentials, which
satisfies simple physical requirements - the potentials are
bounded from below and the weak energy condition is satisfied - as
well as a few more technical assumptions to be discussed below.

\section{Expanding and collapsing models}\label{sec:cosmo}

We consider homogeneous, spatially flat  spacetimes
\begin{equation}\label{eq:g}
g=-\mathrm dt\otimes\mathrm dt+a^2(t)\left[\mathrm dx^1\otimes\mathrm
dx^1+\mathrm dx^2\otimes\mathrm dx^2+\mathrm dx^3\otimes\mathrm dx^3\right]
\end{equation}
where gravity is coupled to a scalar field $\phi$,
self-interacting with a scalar potential $V(\phi)$. Assigning this
function is equivalent to assign the physical properties of the
matter source, and thus the choice of the potential can be seen as
a sort of choice of the equation of state for the matter model.
Therefore, as occurs in models with "ordinary" matter (e.g.,
perfect fluids), the choice of the "equation of state" must be
restricted by physical considerations on stability, positivity of
energy, and so on. We introduce the following conditions:
\begin{enumerate}
\item\label{itm:cond1} $V(\phi)$ is a $\mathcal C^2$ function bounded from below;
\item\label{itm:cond3} The critical points of $V$ are isolated.
They are either minimum points or nondegenerate maximum points.
\end{enumerate}
In the present paper we are going to assume
that $V$ belongs to the set
\begin{equation}\label{as:ipotesi}
\mathfrak V=\{ V:\R\to\R\,:\,\text{\eqref{itm:cond1} and \eqref{itm:cond3} hold}\}.
\end{equation}

The Einstein field equations are given by the following
\begin{subequations}
\begin{align}
&(G^0_0=8\pi T^0_0):\qquad-\frac{3\dot
a^2}{a^2}=-(\dot\phi^2+2V(\phi)),\label{eq:original-efe1}\\
&(G^1_1=8\pi T^1_1):\qquad-\frac{\dot a^2+2a\ddot
a}{a^2}=(\dot\phi^2-2V(\phi)),\label{eq:original-efe2}
\end{align}
\end{subequations}
and have as a consequence the Bianchi identity
\begin{equation}\label{eq:Bianchi}
T^\mu_{\,\,0;\mu}=-2\dot\phi\left(\ddot\phi+V'(\phi)+3\frac{\dot
a}a\dot\phi\right)=0.
\end{equation}

Denoting the energy density of the scalar field by
\begin{equation}\label{eq:energy}
\epsilon=\dot\phi^2+2V(\phi),
\end{equation}
We will actually consider the following system:
\begin{subequations}
\begin{align}
&\left(\frac{\dot a}a\right)^2={\frac\epsilon 3}.\label{eq:initialefe2}\\
&\ddot\phi+V'(\phi)=-3\frac{\dot a}a\,\dot\phi,\label{eq:initialefe1}
\end{align}
\end{subequations}

\bigskip

Of course, it is mandatory to extract the square root in the first
equation. To do this, we introduce the sign function
$$
\chi (t):=\text{sgn} (\dot a(t)).
$$
Physically, this function obviously specifies if the solution
describes an expanding (respectively, collapsing or re-collapsing)
model at time $t$. Thus our final system is composed by
\begin{subequations}
\begin{align}
&\frac{\dot a}a=\chi\sqrt{\frac\epsilon 3}.\label{eq:efe2}\\
&\ddot\phi+V'(\phi)=-\chi\sqrt
3\sqrt\epsilon\,\dot\phi.\label{eq:efe1}
\end{align}
\end{subequations}
Let us observe that, using \eqref{eq:energy} and \eqref{eq:efe1},
the following relation can be seen to hold:
\begin{equation}\label{eq:3}
\dot\epsilon=-2\chi\sqrt 3\sqrt\epsilon\,\,\dot\phi^2,
\end{equation}
that is $\epsilon(t)$ is monotone in each interval in which $\dot
a$ is.

In what follows, we are going to focus on solutions of the
equations \eqref{eq:initialefe2}--\eqref{eq:initialefe1} -- or, better, \eqref{eq:efe2}--\eqref{eq:efe1}. Actually, it can be
proved \cite{jmp} that, if $a\ne 0$ for all $t \geq 0$ and
$(a,\phi)$ are $\mathcal C^2$ functions that solve
\eqref{eq:original-efe1}--\eqref{eq:original-efe2} (with a
non-constant $\phi(t)$) then they are solutions of
\eqref{eq:initialefe2}--\eqref{eq:initialefe1} as well. To prove
the converse, it suffices to derive \eqref{eq:initialefe2} w.r.t. time and use \eqref{eq:initialefe1} to
deduce that if (in some interval $[a,b]$), $\dot a$ is not
identically zero -- that is, $\epsilon\ne 0$ -- then solutions of
\eqref{eq:initialefe2}--\eqref{eq:initialefe1} are also solutions
of \eqref{eq:original-efe1}--\eqref{eq:original-efe2}. Finally, the unique
case in which $\epsilon(t)=0,\forall t>0$ and the Einstein Field
Equations are satisfied is given by the trivial solution
$a(t)=a_0$, $\phi(t)=\phi_0$, with $V(\phi_0)=0$.

\bigskip

We shall divide our discussion of the properties of the solutions
into the two following sections, the first dealing with the
expanding case and the other with the collapsing case. Of course,
the two cases are not disconnected, because in the expanding case,
where $\epsilon(t)$ is decreasing, there might be the possibility
of reaching a vanishing $\epsilon$ in a finite time $t_0$. If this
happens, the model will be ruled from $t_0$ onward by the
equations for the collapsing situation, where
$\epsilon(t)$ is increasing, and it will actually be
proved to be divergent at some finite time for almost every choice
of the initial data, thus yielding a singularity. Moreover this fact, and of course the particular form of system \eqref{eq:efe2}--\eqref{eq:efe1} allows us to perform the same analysis to study also \emph{backwards} qualitative behavior of the solution. Indeed, if the model is expanding in a right neighborhood of $t=0$, and we want to study what happens for $t<0$, replacing $t$ with $-t$ amounts to consider a collapsing model for positive $t$. Then we can conclude that the expanding model ``collapses in the past'', i.e. it originates from a big--bang singularity.

\bigskip

We conclude this section proving a result that will be used
throughout all the paper. It essentially states that the velocity
of  all (finite energy)  solutions which extend infinitely in the
future asymptotically vanishes:
\begin{lem}\label{lem:dotphi0}
Let $\phi(t)$ a solution of \eqref{eq:efe1}
that can be extended for all $t>0$, and such that $\epsilon(t)$ and
$V'(\phi(t))$ are bounded. Then
\begin{equation}\label{eq:limabs}
\lim_{t\to+\infty}\dot\phi(t)=0.
\end{equation}
If, in addition $V''(\phi(t))$ is also bounded, then
\begin{equation}\label{eq:critvalue}
\lim_{t\to +\infty}V'(\phi(t))=0.
\end{equation}
\end{lem}

\begin{proof}
Using \eqref{eq:3} we have
\[
\int_0^{+\infty}\dot\phi(t)^2\,\mathrm dt=\left\vert\int_0^\infty\frac{\dot\epsilon(t)}{2\sqrt 3\sqrt{\epsilon(t)}}\,\mathrm dt\right\vert<+\infty,
\]
so there exists a sequence $t_k\to+\infty$ such that
$\dot\phi(t_k)\to 0$. By contradiction, suppose that
$\exists\rho>0$, and a sequence $s_k\to t_k$ -- that can be taken
such that $s_k>t_k$ -- with $\vert\dot\phi(s_k)\vert>\rho$. Let
$\bar k$ such that $\vert\dot\phi(t_k)\vert<\tfrac\rho 2,\,\forall
k\ge \bar k$, and let $\tau_k,\sigma_k$ sequences such that
$t_k\le \tau_k<\sigma_k\le s_k$, and
\[
|\dot\phi(\tau_k)|=\frac\rho 2,\qquad |\dot\phi(\sigma_k)|=\rho,\qquad\frac\rho 2\le\dot\phi([\tau_k,\sigma_k])|\le\rho.
\]
Note that $\epsilon$ bounded implies $\dot \phi$ bounded, while by assumption, $V'(\phi(t))$ is bounded. Therefore by \eqref{eq:initialefe1},
there exists $M > 0$ such that
\[ |\ddot\phi(t)| \leq M \qquad \forall t.\]
Therefore
\[
\frac\rho 2\leq|\dot\phi(\sigma_k)-\dot\phi(\tau_k)|\leq\int_{\tau_k}^{\sigma_k}|\ddot\phi(t)|\,\mathrm dt\le M(\sigma_k-\tau_k),
\]
that is $\sigma_k-\tau_k\ge\tfrac\rho{2M}$, and therefore
\[
\int_{0}^{+\infty}\dot\phi(t)^2\,\mathrm dt\ge\sum_k\int_{\tau_k}^{\sigma_k}\dot\phi(t)^2\,\mathrm dt\ge\sum_k(\sigma_k-\tau_k)\cdot\left(\frac\rho 2\right)^2\ge\sum_k\frac{\rho^3}{8M},
\]
that diverges. This is a contradiction, and then \eqref{eq:limabs}
must hold.

To prove \eqref{eq:critvalue},
let us first observe that $\exists t_k\to+\infty$ such that $V'(\phi(t_k))\to 0$ -- otherwise,
there would exists $\kappa>0$ such that $\vert V'(\phi(t))\vert\ge\kappa$ definitely, which would imply,
in view of
\eqref{eq:initialefe1}, that $\vert\ddot\phi(t)\vert\ge\kappa/2$, that is absurd since  $\lim_{t\to+\infty}\dot\phi(t)=0$.

Then, let us suppose by contradiction the existence of a constant $\rho>0$,
and a sequence $s_k\to+\infty$ such that $t_k<s_k$ and $\vert V'(\phi(s_k))\vert\ge\rho$.
Therefore, one can choose $\sigma_k,\tau_k$ sequences such that $t_k\le\tau_k<\sigma_k\le s_k$, and
\[
\vert V'(\phi(\tau_k))\vert=\frac\rho 2,\qquad
\vert V'(\phi(\sigma_k))\vert=\rho,\qquad
\frac\rho 2\le\vert V'(\phi(t))\vert\le \rho,\,\forall t\in [\tau_k,\sigma_k].
\]
Then, since by assumption $V''(\phi(t))$ is bounded, there exists
a constant $L>0$ such that
\[
\frac\rho 2=\vert V'(\phi(\sigma_k))-V'(\phi(\tau_k))\vert\le\int_{\tau_k}^{\sigma_k}\vert V''(\phi(t))\vert\,\vert\dot\phi(t)\vert\,\text dt\le L (\sigma_k-\tau_k).
\]
But for sufficiently large $k$ let us observe that
\eqref{eq:initialefe1} implies $\vert\ddot\phi(t)\vert\ge\rho/4$, $\forall t\in [\tau_k,\sigma_k]$, and therefore
$\vert\dot\phi(\sigma_k)-\dot\phi(\tau_k)\vert=\vert\int_{\tau_k}^{\sigma_k}\ddot\phi(t)\text dt\vert\ge\frac{\rho^2}{8L}$
that is a contradiction since $\dot\phi\to 0$.
\end{proof}

\section{Qualitative analysis of the expanding models}\label{sec:espansione}

In this section we study the global behavior of the solutions of
\eqref{eq:efe2}, \eqref{eq:efe1} in the expanding  ($\chi=1$)
case. Regarding the potential, throughout the section we shall
assume that the function $V(\phi)$ belongs to the following subset $\mathfrak E\subseteq\mathfrak V$ \eqref{as:ipotesi}:
\begin{equation}\label{as:ipotesi0}
\mathfrak E=\{V\in\mathfrak V\,:\,\lim_{\vert\phi\vert\to\infty} V(\phi)=+\infty\}
\end{equation}

\begin{rem}\label{rem:plateau}
The above growth condition on the potential is introduced here to avoid cases when $V(\phi)$ admits a flat plateau at infinity; on such situations, the argument of Lemma \ref{valorienergia} below applies, but the function $\phi(t)$ may possibly also diverge as $t\to\infty$, with $\epsilon(t)$ approximating twice the plateau of $V(x)$. In this case, the critical point at infinity is non hyperbolic and arguments exploited, for instance, in Lemma \ref{massimolocale} below do not apply anymore.

With an argument that uses of the Center Manifold Theorem, a result of local asymptotical instability for critical point at infinity is found in \cite{mir2,mir3}, where the potential is supposed to be $V_\infty(1-e^{-\sqrt{3/2}\phi})^2$.
It can be seen that the Center Manifold Theorem applies to extend results of \cite{mir2,mir3} to a more general situation where, as $\phi\to+\infty$, $V(\phi)$ has a finite limit, $V'(\phi)\to 0^+$, and
$$\lim_{\phi\to+\infty} \frac{V''(\phi)}{V'(\phi)}$$
exists finite (and negative).
\end{rem}

We begin our study proving the following:
\begin{lem}\label{valorienergia}
Let $\phi$ be a solution of \eqref{eq:efe1} with $\chi=1$ such that
$\phi(0) = \phi_0$, $\dot \phi(0) =v_0$ with $v_0^2 + 2V(\phi(0))
\equiv \epsilon_0 > 0$ and $(\phi(0), v_0) \neq (\phi_*,0)$ for any
$\phi_*$ critical point of $V$. Then, either there $\exists T =
T(\phi_0,v_0) > 0$ such that $\epsilon(T) = 0$, or $\phi(t)$ is
defined in $\R^+$ and there exists $\phi_*$ critical point of $\,V$
such that:
\[
\lim_{t \to +\infty} \phi(t) = \phi_*, \, \, \lim_{t \to +\infty} \dot \phi(t) = 0.
\]
\end{lem}

\begin{proof}
Due to \eqref{eq:3}, there is $T>0$ such that $\sqrt{3}
\int_0^T{\dot \phi (s)}^2\,\mathrm ds = \sqrt{\epsilon(0)}$ or
\[\int_0^{t}{\dot \phi (s)}^2 \,\mathrm ds < \sqrt{\frac{\epsilon(0)}{3}},\qquad\forall t>0.\]
In the latter case, $\epsilon(t)\in]0,\epsilon_0]$, $\forall t>0$. Since $V$ is bounded from below we obtain that the solution is defined for all $t>0$. Then, by \eqref{as:ipotesi0} we have that $\phi(t)$ is bounded. Therefore, by Lemma \ref{lem:dotphi0},
\[\lim_{t \to +\infty} \dot \phi(t) = 0, \,  \,
\lim_{t \to +\infty} V'(\phi(t)) = 0,\]
and by \eqref{as:ipotesi}, we immediately see that $\phi$ must be
convergent to a critical point of $V$.
\end{proof}

\bigskip

We now proceed to study the two possible cases described in Lemma
\ref{valorienergia}, starting from the behavior of the solutions
for which $\epsilon$ vanishes in a finite time $T>0$. Observe
first that this situation may happen only if $V(\phi(T))\le 0$.
However, we show in the Appendix \ref{sec:eps0} (see remark
\ref{rem:eps0rev}) that if $\inf V<0$ this situation is generic, in the sense
that solutions of this kind exist, and that the set of initial data
leading to solutions such that $V(\phi(T))<0$ is open.
As said before, for $t>T$
the system will be ruled by the equations for the collapsing
situation, discussed in Section \ref{sec:collapse}.

Now, here remains to see what happens if $\phi(t)$ converges to a critical point for $V(\phi)$.
In the following, we shall study first the behavior of the solution near a local
minimum point of $V$.

\begin{lem}\label{minimilocali}
Let $\phi_*$ be a local minimum point of $V$ (not necessarily
nondegenerate).
Assume that there exists a bounded open interval
$]\alpha,\beta[ \ni \phi_*$ such that
\begin{equation}\label{eq:min1}
V'(\phi) > 0,\,\, \forall \phi \in ]\phi_*,\beta], \text{\ and\ }\, V'(\phi) < 0 , \,\,\forall \phi \in [\alpha,\phi_*[;
\end{equation}
\begin{equation}\label{eq:min2}
V(\alpha) = V(\beta) = E > 0;
\end{equation}
\begin{equation}\label{eq:min3}
V(\phi_*) \geq 0.
\end{equation}
Let $\phi$ be a solution such that $0 < \epsilon(0) \leq E$ and
$\phi_0\in[\alpha,\beta]$. Then $\phi$ is defined on all $\R^+$ and
\[
\lim_{t \to +\infty} \phi(t) = \phi_*, \, \lim_{t \to +\infty} \dot \phi(t) = 0.
\]
\end{lem}
\begin{proof}
Let us consider initial data such that $\phi(0) = \phi_0$,
$2V(\phi_0) = E$, $\dot \phi(0) = 0$. Since $\epsilon(t)$ is
non-increasing we have
\begin{equation}\label{eq:sottoE}
\epsilon(t) \leq E, \, \forall t \geq 0.
\end{equation}
Moreover
\[
\ddot \phi(0) = -V'(\phi_0).
\]
If $\phi_0 > \phi_*$, it is $V'(\phi_0) > 0$, then $\phi$ moves
towards $\phi_*$. Analogously if $\phi_0 < \phi_*$, the solution
$\phi$ moves towards $\phi_*$. Now consider an initial data
$(\phi_0, \dot \phi(0))$ such that $\epsilon(0) = E$ and $\dot
\phi(0) \neq 0$, so that $2V(\phi_0) = E-\dot\phi(0)^2 < E$. Note
that, by \eqref{eq:sottoE}, it must be $2V(\phi(t)) \leq E$.
Moreover, if there exists a $\bar t$ such that $2V(\phi(\bar t) =
E$ it must be $\dot \phi(\bar t) = 0$ and arguing as above we see
that $\phi$ moves towards $\phi_*$. Therefore $\alpha \leq \phi(t)
\leq \beta$ for all $t$. Then, by Lemma \ref{valorienergia} we
deduce that $\lim\limits_{t \to +\infty}\phi(t) = \phi_*$ which is
the only critical point of $V$ in the interval $[\alpha,\beta]$,
and $\lim\limits_{t \to +\infty}\dot\phi(t) = 0$.
\end{proof}

\begin{rem}\label{oscillazione}
Whenever $V(\phi_*) = 0$ further information about the asymptotic
behavior of the solution near the minimum point $\phi_*$ can be
obtained \cite{r2}. In that paper it is indeed proved the
oscillation of $\phi$ around the minimum point using asymptotic
analysis. If $\phi_*$ is a nondegenerate minimum point (which is
the more relevant physical case) we can prove the oscillatory
character also using an alternative argument as follows. Suppose,
by contradiction, that $\phi$ does not oscillate around $\phi$,
this means that there exists $t_1 > 0$ such that
\[
\phi(t) \neq \phi_*, \, \forall t \geq t_1.
\]
To fix our ideas suppose $\phi(t) > \phi_* \, \forall t\geq t_1$. Then, if $t \geq t_1$ and $\dot \phi(t) = 0$ it is $\ddot\phi(t) = -V'(\phi(t)) < 0$ (since
$V' > 0$ in a right neighborhood of $\phi_*$). This implies that there exists $t_2 > 0$ such that
\[
\dot\phi(t) < 0, \, \forall t \geq t_2 .
\]
Moreover there are not sequences $t_k \rightarrow +\infty$ where
$\ddot \phi(t_k) = 0$. Otherwise we should have $\dddot \phi(t_k)
= -(V''(\phi(t_k) + \sqrt{3}\frac{\dot
\epsilon}{\sqrt{\epsilon}})\dot \phi(t_k) = -(V''(\phi(t_k) -
12(\dot \phi(t_k))^2)\dot \phi(t_k))$, so $\dddot \phi(t_k) < 0$
for any $k$ sufficiently large, getting a contradiction. Therefore
there exists $t_3 > 0$ such that
\[
\ddot\phi(t) > 0, \, \forall t \geq t_3.
\]
Then, by \eqref{eq:efe1},
\[
-V'(\phi) - \sqrt{3\epsilon}{\dot \phi} > 0, \, \forall t \geq t_3,
\]
that implies (recalling that $\dot \phi < 0$)
\[
3{(\dot \phi)}^4 +6V(\dot\phi)^2 -(V')^2 \geq 0, \, \forall t \geq t_3
\]
and
\[
(\dot\phi)^2 \geq \sqrt{V^2 + \frac{(V')^2}3} - V, \, \forall t \geq t_3.
\]
Now $V(\phi) = a(\phi-\phi_*)^2 +o(\phi - \phi_*)^2$ with $a>0$, so there exists $b > 0$ and $t_4 > 0$ such that
\[
|\dot \phi| \geq b\sqrt{\phi-\phi_*}, \, \forall t_4,
\]
(recall that $\phi(t) > \phi_*$). Finally $\forall t \geq t_4$
\[
\dot \phi \leq -b\sqrt{\phi-\phi_*},
\]
and integrating we obtain,
\[
\sqrt{\phi(t) -\phi_*} - \sqrt{\phi(t_4) -\phi_*} \leq -\frac{b}{2}t, \forall t\geq t_4,
\]
in contradiction with $\phi(t) \geq \phi_*$ for any $t \geq t_4$.
\end{rem}

Now we shall study the behavior near a local nondegenerate maximum
point of $V$.
\begin{lem}\label{massimolocale}
Let $\phi_*$ be a local nondegenerate maximum point of $V$ such
that $V(\phi_*) \geq 0$. Then there are only two solutions
$\phi_1$ and $\phi_2$ of \eqref{eq:efe1} such that $\lim\limits_{t
\rightarrow +\infty}\phi_i(t) = \phi_*$, $\lim\limits_{t
\rightarrow +\infty}\dot \phi_i(t) = 0$, $i=1,2$. Moreover there
exists $t_* > 0$ such that $\phi_1(t) > \phi_*$ for all $t \geq
t_*$ and $\phi_2(t) < \phi_*$ for all $t \geq t_*$.
\end{lem}
\begin{proof}
Consider the first order system associated to \eqref{eq:efe1}
\begin{equation}\label{sistemanonlineare}
\begin{cases}
& \dot x= -V'(y) - \sqrt{E(x,y)}x \\
& \dot y = x,
\end{cases}
\end{equation}
where $E(x,y) = x^2 + 2V(y)$.

Let $(0,y_*)$ be an equilibrium point of \eqref{sistemanonlineare}, namely $V'(y_*)=0$. Suppose
\[
E(0,y_*) = 2V(y_*) > 0.
\]
The linearized system of \eqref{sistemanonlineare} at $(0,y_*)$ is
\begin{equation}\label{sistemalineare}
\begin{cases}
& \dot \alpha= -V''(y_*)\beta - \sqrt{2V(y_*)}\alpha \\
& \dot \beta = \alpha.
\end{cases}
\end{equation}

Since, by our assumption $V''(y_*) < 0$, the matrix
$ \left(
\begin{matrix}
-\sqrt{2V(y_*)} & -V''(y_*)\cr
1&0\cr
\end{matrix}
\right)$
has one positive eigenvalue and one negative eigenvalue. Therefore the proof follows immediately by classical theory on stability of nonlinear dynamical
system in a neighborhood of an equilibrium point (cf. e.g. \cite{perko}).

Note that the function $\sqrt{E(x,y)}x$ is not of class $C^1$ in a neighborhood of $(0,y_*)$ whenever $V(y_*) = 0$. However the quantity $\sqrt{E(x,y)}x$
is an infinitesimal of order greater than 1 in $(0,y_*)$, therefore the classical theory can be adapted to this case, studying the linear system
\[
\begin{cases}
& \dot \alpha= -V''(y_*)\beta  \\
& \dot \beta = \alpha,
\end{cases}
\]
obtaining the same conclusion of the case $V(y_*) > 0$.
\end{proof}

\bigskip

Collecting the results shown in this Section (and in Appendix \ref{sec:eps0}) we obtain the following
\begin{prop}
Suppose that $V$ belongs to the set $\mathfrak E$ \eqref{as:ipotesi0}. Let be $\phi$ a solution of \eqref{eq:efe1} where $\chi=1$, with
initial data $\phi(0) = \phi_0$, $\dot \phi(0) = v_0$, and $v_0^2
+ 2V(\phi_0) > 0.$ Then, one of the two mutually exclusive
situations occur:

\smallskip

\noindent either
\begin{equation}\label{eq:zero}
\text{$\exists T=T(\phi_0,v_0) > 0$ such that $\epsilon(T)= (\dot \phi(T))^2 +2V(\phi(T)) = 0$,}
\end{equation}
and the set of initial data such that \eqref{eq:zero}
is satisfied with $V(\phi(T)) < 0$ is open and not empty,
while there are only a finite number of solutions satisfying \eqref{eq:zero} with
$V(\phi(T)) = 0$,

\smallskip
\noindent or
\begin{equation}\label{eq:critico}
\text{$\exists\phi_*$ critical point of $V$ with $V(\phi_*)\geq 0$ such that $\lim\limits_{t \rightarrow +\infty}\phi(t) = \phi_*$,
and $\lim\limits_{t \rightarrow +\infty}\dot \phi(t) = 0$,}
\end{equation}
and the measure of the initial data set such that $\phi(t) \rightarrow \phi_*$ maximum point is zero, while the set of initial data such that
$\phi(t) \rightarrow \phi_*$ minimum point is open and not empty.

Moreover, situation \eqref{eq:zero} occurs only if $\inf V<0$.
\end{prop}

\section{Qualitative analysis of the collapsing models
}\label{sec:collapse}

The aim of the present section is to study the qualitative
behavior of the solution of \eqref{eq:efe2}--\eqref{eq:efe1} in
the collapsing case ($\chi=-1$), for which we will require that
$V(\phi)$ satisfies some further conditions. We first
assume that
$\exists V^*$ positive constant such that:
\begin{align}
&\text{the set} \quad B:=\{\phi\in\R\,:\,V(\phi)\le V^*\} \quad \text{is bounded,} \label{eq:bound} \\
\intertext{and}
&\phi\ge\sup B\Longrightarrow V'(\phi)>0,\qquad
\phi\le\inf B\Longrightarrow V'(\phi)<0.\label{eq:Vprime}
\end{align}
Moreover, defining the function
\begin{equation}\label{eq:Y}
u(\phi):=\frac{V'(\phi)}{2\sqrt 3V(\phi)},
\end{equation}
and introducing the growth conditions
\begin{align}
&\limsup_{\phi\to\pm\infty} |u(\phi)|< 1, \label{eq:u-infty}\\
&\exists\lim_{\phi\to\pm\infty} {u'(\phi)} \, (=0).\label{eq:uprime-infty}
\end{align}
we assume that $V$ belongs to the subset $\mathfrak C\subseteq\mathfrak V$ defined as
\begin{equation}\label{as:assum2}
\mathfrak C=\{V\in\mathfrak V\,:\, V \text{\ satisfies \eqref{eq:bound}--\eqref{eq:Vprime} and \eqref{eq:u-infty}--\eqref{eq:uprime-infty}}\}.
\end{equation}
We stress that this class contains the potentials which are
usually considered in the physical applications. For instance, it
contains the standard quartic potential ($V(\phi)=-\frac12
m^2\phi^2+\lambda^2\phi^4$) and, more generally, all the
potentials bounded from below whose asymptotic behavior is
polynomial (i.e. $\phi^{2n}$ for $|\phi|\to\infty$, $n \geq1$). Moreover, it obviously suffices for $V\in\mathfrak C$ to diverge at infinity to be also in the class $\mathfrak E$ studied in previous Section \ref{sec:espansione}.

Let $t=0$ be the initial time and let $\phi_0=\phi(0)
\in B$. We consider only data which satisfy the weak energy
condition and thus, we assume
\begin{equation}\label{eq:wec}
\epsilon_0=\epsilon(0)=\dot\phi^2(0)+2V(\phi(0))\ge 0.
\end{equation}
We start considering initial data such that $\epsilon_0$ is
"sufficiently large", more precisely
\begin{equation}\label{eq:initial}
\epsilon_0\ge 2V^*,
\end{equation}
where $V^*$ has been defined just above, and proceed to show that the solution $\phi(t)$ diverges in a
finite time for almost every choice of the initial data
(Proposition \ref{thm:sing}). Then, we will  extend the result by
removing the restriction \eqref{eq:initial}.

 Before starting the proof, let us notice the following:

\begin{rem}\label{rem:dotphi0}
If $\dot\phi(0)=0$ then $\dot\phi(t)\ne 0$ in a right neighborhood of $t=0$. Indeed, $2V(\phi_0)=\epsilon_0\ge 2V^*$, then $\phi_0\not\in\overset\circ B$, the set of internal points of $B$. This implies by \eqref{eq:Vprime} that $V'(\phi_0)\ne 0$, so that \eqref{eq:efe1} shows $\ddot\phi(0)\ne 0$.
\end{rem}

\begin{rem}\label{rem:doteps}
Recall from \eqref{eq:3} that $\epsilon(t)$ is monotonically non decreasing. Actually, we can also observe that, $\forall t>0$, $\epsilon(t)$ is strictly greater that $\epsilon_0$: it is straightforward if $\dot\phi(0)>0$, but is also true, in view of Remark \ref{rem:dotphi0}, if $\dot\phi(0)=0$.
\end{rem}

\begin{rem}\label{rem:delta}
Given $t_0>0$, there exists a constant $\delta>0$ such that, for all $t_1,t_2$ such that $t_0<t_1<t_2$, and $\phi([t_1,t_2])\subseteq\overset\circ B$,
$\dot\phi(t)^2\ge\delta,\,\forall t\in [t_1,t_2]$. Indeed, let $\delta=\epsilon(t_0)-\epsilon_0>0$
(cf. Remark \ref{rem:doteps}). Then, if $t\in [t_1,t_2]$, $\phi(t)^2=\epsilon(t)-2V(\phi(t))\ge\epsilon(t_1)-2V(\phi(t))\ge\delta+\epsilon_0-2V^*\ge\delta$,
where we have used \eqref{eq:bound} and \eqref{eq:initial}.
\end{rem}

\begin{rem}\label{rem:maxloc}
Let $t>0$ such that $\phi(t)\not\in\overset\circ B$ and $\dot\phi=0$. Then, \eqref{eq:efe1} and \eqref{eq:Vprime} imply that $t$ is a local extremum for $\phi$. In particular, it is a local maximum if $\phi(t)>\sup B$, and is a local minimum if $\phi(t)<\inf B$.
\end{rem}

Let $\I\subseteq[0,+\infty)$ be the maximal right
neighborhood of
$t=0$ where the solution $\phi(t)$ is defined, and call
\[
t_s=\sup\I\in\R^+\cup\{+\infty\}.
\]
We can now state the following lemma.

\begin{lem}\label{thm:phidotphi}
If $\phi(t)$ is bounded, also $\dot\phi(t)$ is bounded, and
$t_s=+\infty$.
\end{lem}
\begin{proof}
Let $K$ a bounded set such that $\phi(t)\in K,\,\forall t\in \I$, and
argue by contradiction, assuming $\dot\phi$ not bounded, and supposing that $\sup_\I\dot\phi(t)=+\infty$ -- the same argument can be used if $\dot\phi$ is  unbounded only below. Let $\bar t$ such that
\[
\dot\phi(\bar t)>\frac{\sup_{K}V'(\phi)}{\sqrt{3\epsilon_0}},
\]
therefore \eqref{eq:efe1} implies that $\ddot\phi(\bar t)>0$, and then $\dot\phi$ is increasing in $\bar t$. But this means that $\dot\phi$ is increasing
in $\bar t$, therefore \eqref{eq:efe1} shows that $\ddot\phi\ge 0,\,\forall t\geq\bar t$, and then $\dot\phi$ is eventually increasing, namely $\lim_{t\to t_s}\dot\phi(t)=+\infty$.
This shows that $t_s<+\infty$, otherwise it would be $\phi(t)-\phi_0=\int_0^t\dot\phi(\tau)\,\mathrm d\tau$, which would diverge as $t\to+\infty$. Then $t_s\in\R$.

With the position
\[
\lambda(t)=\frac1{\dot\phi(t)},
\]
equation \eqref{eq:efe1} implies
\[
\dot\lambda(t)=\frac{V'(\phi(t))}{\dot\phi(t)^2}-\sqrt 3\sqrt{1+\frac{2V(\phi(t))}{\dot\phi(t)^2}}\xrightarrow{t\to t_s^-}-\sqrt 3,
\]
and then $\lambda(t)$, near $t=t_s$, behaves like $\sqrt 3(t_s-t)$. Since $\dot\phi(t)$ positively diverges,  $\lim_{t\to t_s}\phi(t)=\phi^*$ exists, and is finite since $\phi(t)$ is bounded by hypothesis. But the quantity
\[
\phi(t)-\phi^*=\int_{t_s}^t\dot\phi(\tau)\,\mathrm d\tau=\int_{t_s}^t\frac1{\lambda(\tau)}\,\mathrm d\tau
\]
diverges, which is absurd. This shows that $\dot\phi(t)$ must be bounded too.

To show that $t_s=+\infty$, we proceed again by contradiction. Let
$t_k\to t_s$ be a sequence such that the sequence
$(\phi(t_k),\dot\phi(t_k))$ converges to a finite limit
$(\phi_*,\dot\phi_*)$. Solving the Cauchy problem \eqref{eq:efe1}
with initial data $(\phi_*,\dot\phi_*)$ shows that the solution
$\phi(t)$ is $\mathcal C^1$, the solution could be prolonged on a
right neighborhood of $t_s$.
\end{proof}

\begin{prop}\label{thm:limsup}
The function $\phi(t)$ is unbounded.
\end{prop}
\begin{proof}
Let by contradiction $\phi(t)$ be bounded. Then, by Lemma \ref{thm:phidotphi}, $\vert\dot\phi(t)\vert\le M$ for a suitable constant $M$,
and $t_s=\sup\I=+\infty$. In particular $V(\phi(t)), V'(\phi(t)), \epsilon(t),\ddot\phi(t)$ are bounded also.
Then Lemma
\ref{lem:dotphi0}
 says that $\lim_{t\to+\infty}\dot\phi(t)=0$, that implies
$\lim_{t\to+\infty}2V(\phi(t))=\lim_{t\to+\infty}\epsilon(t)>\epsilon_0\ge 2V_*$, and so $\phi(t)\not\in B$
for any $t$ sufficiently large. Then, for large $t$, $\phi(t)$ moves in a region where $V$ is invertible. Since $V(\phi(t))$ converges this happen also
for $\phi(t)$ which converges to a point $\phi_* \not \in \overset\circ B$. Then by
\eqref{eq:efe1}, $\ddot\phi(t)\to-V'(\phi_*)\ne 0$ (since $\dot \phi (t) \rightarrow 0$), that is absurd. Then $\limsup_{t\to t_s}\vert\phi(t)\vert=+\infty$.
\end{proof}

We have shown that $\phi(t)$ is not bounded.
Now we want to show that, actually, $\phi(t)$ monotonically diverges in the approach to $t_s$.

\begin{rem}\label{rem:rho}
We observe that the quantity
\begin{equation}\label{eq:rho}
\rho(t):=\frac{2V(\phi(t))}{\dot\phi(t)^2},
\end{equation}
satisfies the equation
\begin{equation}\label{eq:ODErho}
\dot\rho=2\sqrt3\dot\phi\,\rho\sqrt{(1+\rho)}\left(u(\phi)\sqrt{1+\rho}-\mathrm{sgn}\,(\dot\phi)\right).
\end{equation}
\end{rem}

\begin{prop}\label{thm:div}
Let $V\in\mathfrak C$ \eqref{as:assum2}, and
let $\phi:\I\subseteq\R\to\R$ the solution of \eqref{eq:efe1}. Then $\dot\phi(t)$ is eventually not zero, and
\[
\lim_{t\to\sup\I}|\phi(t)|=+\infty.
\]
\end{prop}
\begin{proof}
By contradiction, let $s_k$ be a sequence of local minima, and $t_k$ a sequence of local maxima for $\phi(t)$,
both convergent to $t_s$, and such that $s_k<t_k<s_{k+1}<t_{k+1}$, for each $k$ (recall that $|\phi(t)|$ is unbounded). Then Remark \ref{rem:doteps} implies
\[
V_*\le\frac{\epsilon_0}2\le\frac12\epsilon(s_k)\le\frac12\epsilon(t_k)\le
\frac12\epsilon(s_{k+1})\le\frac12\epsilon(t_{k+1}),
\]
that is
\[
V_*\le V(s_k)\le V(t_k)\le V(s_{k+1})\le V(t_{k+1}).
\]
Recalling Remark \ref{rem:maxloc}, $s_k\le\inf B$ and $t_k\ge \sup
B$, that is $\phi(t)$ crosses $B$ infinitely many times. We claim
that $\lim_{t\to t_s}\epsilon(t)=+\infty$. Otherwise,
$\epsilon(t)$ would be bounded, and also $\dot\phi(t)$ would be.
In particular, taken $\sigma_k,\tau_k$ sequences converging to
$t_s$ such that $\sigma_k<\tau_k$, and
\[
\phi(\sigma_k)=\inf B,\qquad\phi(\tau_k)=\sup B,\qquad\phi([\sigma_k\tau_k])\subseteq B,
\]
and called
$
N=\sup B-\inf B,
$
it must be
\[
N=\phi(\tau_k)-\phi(\sigma_k)=\int_{\sigma_k}^{\tau_k}\dot\phi(t)\mathrm dt\le M(\tau_k-\sigma_k),
\]
and so $\tau_k-\sigma_k\ge \tfrac NM$. Then, by \eqref{eq:3}, and Remark \ref{rem:delta},
\[
\sqrt{\epsilon(\tau_k)}-\sqrt{\epsilon(\sigma_k)}=\int_{\sigma_k}^{\tau_k}
\frac{\dot\epsilon}{2\sqrt\epsilon}\,\mathrm dt=
\sqrt 3\int_{\sigma_k}^{\tau_k}\dot\phi^2\,\mathrm dt\ge\delta\sqrt {3}(\tau_k-\sigma_k)\ge\delta\sqrt 3\frac NM,
\]
and then $\epsilon(t)$ must diverge.

Now, recalling \ref{eq:u-infty}, let $\theta>0$ be a sufficiently
small constant, and let $r_k\to t_s$ a sequence such that
$\phi(r_k)=\bar\phi$, with $u(\phi) \leq 1-\theta$ as
$\phi\ge\bar\phi$. Then
$\dot\phi(r_k)^2=\epsilon(r_k)-2V(\bar\phi)\to+\infty$. Let us
suppose that $\dot\phi(r_k)$ is unbounded above -- analogously one
can argue if it is only unbounded below. Up to subsequences, we
can suppose that $\dot\phi(r_k)$ positively diverges. Moreover,
recalling \eqref{eq:rho}, $\rho(r_k)\to 0$ and then, if $k$ is
sufficiently large,
\[
\vert\rho(r_k)\vert<\left(\frac1{1-\theta}\right)^2-1,
\]
so that, since $\phi(r_k)=\bar\phi$ and $u(\bar\phi) \leq 1-\theta$, we have
\[
u(\phi(r_k))\sqrt{1+\rho(r_k)}-1<0.
\]
Then, $\dot\rho(r_k)<0$, that is $\rho(t)$ is decreasing at $r_k$. But we can observe that $\phi(t)$ is increasing in $r_k$, ensuring $u(\phi(t))<1-\theta$, and so, in a right neighborhood of $r_k$,
\[
u(\phi(t))\sqrt{1+\rho(t)}-1<(1-\theta)\sqrt{1+\rho(t)}-1
\]
that is decreasing. We conclude that the function $\rho(t)$ decreases for $t>r_k$, until $t$ equals a local maximum $t_k$, where $\dot\phi$ vanishes.
But since $2V(\phi(t))\le\kappa\dot\phi(t)^2$ for some $\kappa>0$, this fact would imply $\epsilon(\phi(t_k))=0$, a contradiction.
Then $\dot \phi$ is eventually non zero. This implies that there exists $\lim\limits_{t\rightarrow +\infty}|\phi(t)|$ so,
by Proposition \ref{thm:limsup}, it is $+\infty$.
\end{proof}

We have shown so far that $\phi(t)$ diverges, as $t$ approaches
$t_s=\sup\I$. In the following, we will show that $t_s\in\R$
for almost every solution, in the sense that there exists a set of
initial data, dense in the set of the admissible ones, such that
the solution is defined until a certain finite comoving time $t_s$
(depending on the data). Henceforth we will suppose (just to fix our ideas) $\phi(t)$
\emph{positively} diverging and $\dot \phi(t) > 0$  $\forall t\ge\bar
t$. Then, in the interval $[\bar t,+\infty)$, $\rho$ can be seen
as a function of $\phi$, that satisfies by \eqref{eq:ODErho} the
ODE
\begin{equation}\label{eq:eva}
\frac{\mathrm d\rho}{\mathrm d\phi}=2\rho\sqrt{3(1+\rho)}\left(u(\phi)\sqrt{1+\rho}-1\right).
\end{equation}

\begin{rem}\label{rem:negative}
It is not restrictive to study \eqref{eq:eva} for large and positive $\phi$.
Indeed, consider the case $\phi\to-\infty$, and take $\widetilde\rho(\psi)=
\rho(-\psi)$, $\widetilde V(\psi)=V(-\psi)$. Then, \eqref{eq:ODErho} gives
\begin{multline*}
\frac{\mathrm d\widetilde\rho}{\mathrm d\psi}(\psi)=
-\frac{\mathrm d\rho}{\mathrm d\psi}(-\psi)=-\rho(-\psi)2\sqrt{3(1+\rho(-\psi))}\left(
\frac{V'(-\psi)}{2\sqrt 3V(-\psi)}\sqrt{1+\rho(-\psi)}+1\right)=\\
\widetilde\rho(\psi)2\sqrt{3(1+\widetilde\rho(\psi))}\left(\frac{\widetilde V'(\psi)}{2\sqrt 3\widetilde V(\psi)}\sqrt{1+\widetilde\rho(\psi)}-1\right),
\end{multline*}
where $\psi$ is large and positive. Once we observe that $\widetilde V$ satisfies, for large and positive $\psi$, the same assumptions as $V$, we have that \eqref{eq:eva} also controls the behavior of $\rho$ as $\phi\to-\infty$.
\end{rem}

We now state the following crucial result.

\begin{lem}\label{thm:perla}
Except at most for a measure zero set of initial data, the function $\rho(t)$
goes to zero for $t\to\sup\I$.
\end{lem}
\begin{proof}
The proof will be carried on by studying qualitatively the solutions of the ODE \eqref{eq:eva}. By virtue of Proposition \ref{thm:div}, and Remark \ref{rem:negative}, we will be interested in those solutions which can be indefinitely prolonged on the right.

With the variable change $y=\sqrt{1+\rho}$, \eqref{eq:eva} becomes
\begin{equation}\label{eq:y}
\frac{\mathrm dy}{\mathrm d\phi}={\sqrt 3} (y^2-1)\left(u(\phi)y-1\right),
\end{equation}
where we recall that $u(\phi)$ is given by \eqref{eq:Y}.

Let us consider solutions of \eqref{eq:y} defined in $[\phi_0,+\infty)$.
If $\liminf_{\phi\to+\infty}{\vert y(\phi)}{u(\phi)-1\vert}>0$ then, necessarily,
$u(\phi)y(\phi)-1< 0$ (and $y(\phi)$ decreases, and goes to 1), otherwise the solution would not be defined in a neighborhood of $+\infty$.
In short, if $y(\phi)$ is a solution defined in a neighborhood of $+\infty$, the following behaviors are allowed:
\begin{enumerate}
\item\label{itm:1} either $y(\phi)$ is eventually weakly decreasing, and tends to 1 (so that $\rho$ tends to 0),
\item\label{itm:2} or $\liminf_{\phi\to+\infty}{\vert y(\phi)}{u(\phi)-1\vert}=0$.
\end{enumerate}
With the variable and functions changes
\[
s=e^{-{\sqrt 3}\phi},\qquad z=uy-1,
\]
equation \eqref{eq:y} takes the form
\begin{equation}\label{eq:z}
(s \,u)\dot z(s)=-(z+1)^2 z + (u^2+s \,\dot u) z+ s \,\dot u.
\end{equation}
But it can be easily seen,
using \eqref{eq:uprime-infty}, and the identity
$$s \,\dot u(s)=-\tfrac1{\sqrt 3}u'(\phi),$$ that \eqref{eq:z} satisfies
the assumption of Theorem \ref{thm:ode} in  Appendix \ref{sec:proof}, and so there exists a unique solution for the Cauchy problem given by \eqref{eq:z} with the initial condition $z(0)=0$, that is furthermore the only possible solution of the ODE \eqref{eq:z} with $\liminf_{s\to 0}\vert z(s)\vert=0$, and this fact results in a unique solution satisfying case \eqref{itm:2} above, whereas  all other situations lead to case \eqref{itm:1}.
\end{proof}

Now, we show that the singularity forms in a finite amount of
time, for almost every choice of initial data.

\begin{prop}\label{thm:sing}
If $\epsilon_0\ge 2 V_*$, except at most for a measure zero set of initial data, there exists $t_s<+\infty$ such that $t_s=\sup\I$, and the following facts holds:
\begin{align}
&\lim_{t\to t_s^-}\epsilon(t)=+\infty,\label{eq:eps-div}\\
&\lim_{t\to t_s^-}\dot\phi(t)=+\infty,\label{eq:dot-div}\\
&\lim_{t\to t_s^-}a(t)=0.\label{eq:a-van}
\end{align}
\end{prop}

\begin{proof}
From \eqref{eq:energy} and \eqref{eq:efe2}--\eqref{eq:efe1}, it easily follows that
\begin{equation}\label{eq:dotenergy}
\dot\epsilon(t)=2\sqrt 3\epsilon^{3/2}\frac1{1+\rho(t)}.
\end{equation}
But Lemma above ensures that $\rho$ (defined in \eqref{eq:rho}) goes to zero for almost every choice of initial data.
For this choice, $\dot\epsilon(t)>\epsilon^{3/2}$ in a left neighborhood of $\sup\I$. Then
$t_s < +\infty$ and
\eqref{eq:eps-div} easily follows, using comparison theorems in ODE.

We already know that $\phi(t)\to+\infty$. To prove \eqref{eq:dot-div} we observe
that
\[
\frac{\epsilon(t)}{\dot\phi(t)^2}=\rho(t)+1\to 1,
\]
and then \eqref{eq:eps-div} implies \eqref{eq:dot-div}.

To prove \eqref{eq:a-van}, we first prove that
\begin{equation}\label{eq:inteps-as}
\lim_{t\to t_s^-}\int_0^{t}\sqrt\epsilon\,\mathrm d\tau=+\infty.
\end{equation}
Indeed, recalling that $\rho(t)$ is eventually bounded, $\exists \bar t>0$ such that $\epsilon(t)=\dot\phi^2(t)(1+\rho(t))\ge\kappa^2\dot\phi(t)^2$, for some suitable constant $\kappa>0$, so $\sqrt{\epsilon(t)}>\kappa\dot\phi(t)$, $\forall t>\bar t$, and then
\[
\int_{\bar t}^t\sqrt\epsilon\,\mathrm d\tau\ge\kappa\int_{\bar t}^t \dot\phi(\tau)\,\mathrm
d\tau=\kappa(\phi(t)-\phi(\bar t)).
\]
The righthand side above diverges because of Proposition \ref{thm:div}, and  \eqref{eq:inteps-as} is proved. Moreover, using \eqref{eq:efe2},
\[
\log\frac{a(t)}{a(0)}=\int_0^t\frac{\dot a}a\,\mathrm d\tau=-\frac1{\sqrt
3}\int_0^t\sqrt\epsilon\,\mathrm d\tau,
\]
from which \eqref{eq:a-van} follows, since by \eqref{eq:inteps-as} the righthand side above negatively diverges
as $t\to t_s^-$.

\end{proof}

\subsection{The case $\epsilon_0< 2V^*$}\label{sec:sottosoglia}

Now we are going to consider the case where \eqref{eq:initial} does not hold.
This means that there exists at least one critical point
$\phi_*$ for $V(\phi)$ such that
\begin{equation}\label{eq:initial2}
2V^*>V(\phi_*)\geq\epsilon_0.
\end{equation}
First of all, we take care of the case  $\epsilon_0=0$. Since the
potential can be negative, one should consider the case in which
the energy vanishes "dynamically" (that is, $\phi$ solves the
equation $\dot\phi(t)^2=-2V(\phi(t))\,\forall t$). However, these
functions are not solutions of Einstein field equations
\eqref{eq:original-efe1}--\eqref{eq:original-efe2}, since
$\epsilon(t)\equiv 0$ implies $a(t)=a_0$ and then
$\dot\phi(t)=0,\,V(\phi(t))=V(\phi_0)=0$. Thus, one is left with a set
of constant solutions with zero energy of the form $(a_0,\phi_0)$,
with $a_0$ positive constant and $V(\phi_0)=0$. This shows that at
the "boundary of the weak energy condition" local uniqueness of
the field equations may be violated if $V'(\phi_0)\ne 0$. However the set of the
initial data for the expanding equations intersecting such
solutions in a finite time has zero measure (note that the points
$\phi_0$ such that $V(\phi_0)=0$ and $V'(\phi_0)\neq 0$ are
isolated). Further, we prove (see Lemma \ref{lem:eps0} in Appendix
\ref{sec:eps0}) a result of local existence and uniqueness for the
solutions of equation \eqref{eq:efe1} satisfying $\epsilon(0)=0$,
but for which $\epsilon(t)>0$ for $t>0$.

Let us now study the evolution of $\phi(t)$ assuming $\epsilon(t)<
2V^*$, $\forall t\in\I$ (otherwise the previously obtained results
would apply). In this case $\phi(t) \in B$ for any $t \geq 0$ and,
since $\dot \phi$ is bounded,  $\sup\I=+\infty$. Then  Lemma
\ref{lem:dotphi0} applies, to find $\lim_{t\to
+\infty}\dot\phi(t)=0,\,\lim_{t\to +\infty}V'(\phi(t))=0.$
Moreover by assumption \ref{as:ipotesi} there exists $\phi_*$
critical point of $V$ such that $\lim_{t\to
+\infty}\phi(t)=\phi_*$, (with critical value $V(\phi_*) =
\lim_{t\to +\infty}\epsilon(t) \in ]0,2V_*[$.

\bigskip

Whenever $\phi_*$ is a (nondegenerate) maximum point we can study
the linearization of the first order system equivalent to
\eqref{eq:efe1} in a neighborhood of the equilibrium point
$(0,\phi_*)$, as done in Lemma \ref{massimolocale}, obtaining a
result totally analogous to Lemma \ref{massimolocale}. Moreover,
by the results of Lemma \ref{minimilocali} (which can be seen as
the "time-reversed" version) we see if $\phi_*$ is a minimum point
and $\phi$ starts with initial data close to $(0,\phi_*)$, then
$\phi$ moves far away from $\phi_*$.

Therefore, under the assumptions made, for almost every choice of
the initial data, the function $\phi(t)$ must be such that its
evolution cannot be contained in the compact set $B$ defined in
\eqref{eq:bound}, namely, $\epsilon(t_1)\ge 2V_*$ for some $t_1>0$,
which allows us to apply the theory we already know to show that
the singularity forms for almost every choice of the initial data.
The results are summarized in the following theorem:

\begin{teo}\label{thm:sing-gen}
Suppose that $V(\phi)$ belongs to the set $\mathfrak C$ \eqref{as:assum2}.
Then, except at most for a measure zero set of initial data satisfying weak energy condition, there exists $t_s\in\R$ such that the scalar field
solution becomes singular at $t=t_s$, that is
$\lim_{t\to t_s^-}\epsilon(t)=+\infty,$ and $\lim_{t\to t_s^-}a(t)=0.$ Moreover
$\lim_{t\to t_s^-}\dot\phi(t)=+\infty.$

\end{teo}

\subsection{Examples}
\begin{ex}\label{ex:poly}
The above results hold for all potentials with polynomial
leading term at infinity (i.e. $\lambda^2\phi^{2n}$).
\begin{figure}
\begin{center}
\psfull \epsfig{file=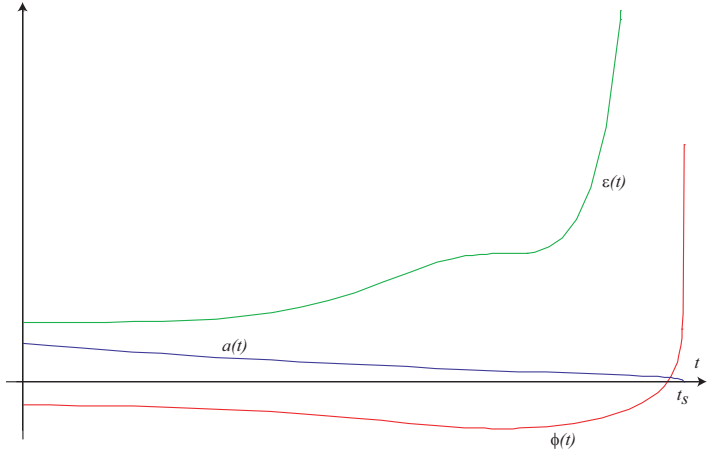, height=6cm}
\caption{Behavior of the function $\phi(t),\epsilon(t)$,
and $a(t)$ with potential given by $V(\phi)=1-\phi^2+\phi^4$.
The initial conditions are $\phi_0=-0.6,\,\dot\phi_0=0,\,a_0=1$. The time of collapse $t_s$ approximately equals $t_s=2.1$.}\label{fig:poly}
\end{center}
\end{figure}
For instance, for a quartic potential $V(\phi)=-\frac12
m^2\phi^2+\lambda^2\phi^4$, with $\lambda,m\ne 0$, the function
$u(\phi)$ goes as $\tfrac 2{\sqrt 3\phi}$ for $\phi\to +\infty$,
and all conditions listed above hold. The behavior for a
particular example from this class of potential is given in Figure
\ref{fig:poly}.
\end{ex}

\begin{ex}\label{ex:exp}
In the case of exponential potentials, the results hold for
asymptotic behaviors with leading term at infinity of the form $V_0
e^{2\sqrt{3} \lambda|\phi|}$ with $\lambda <1$. For instance for
$V(\phi)=V_0 e^{2\sqrt 3 \lambda\sqrt{\phi^2+\gamma^2}}$, where
$V_0,\lambda,\gamma>0$, the quantity $u(\phi)$ goes like
$\lambda$, and so \eqref{eq:u-infty} is verified if $\lambda<1$.
\begin{figure}
\begin{center}
\psfull \epsfig{file=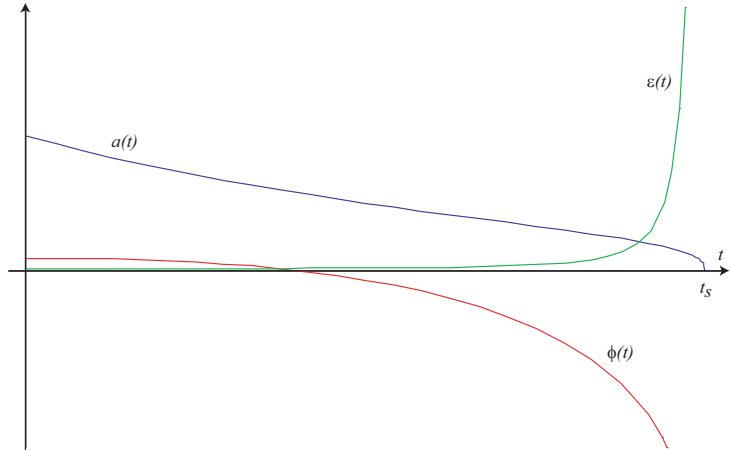, height=6cm}
\caption{Same as Figure \ref{fig:poly}, but now the potential is given by $V(\phi)=e^{2\sqrt{3/2(\phi^2+1)}}/100$. The initial conditions are $\phi_0=0.1,\,\dot\phi_0=0,\,a_0=1$. The time of collapse $t_s$ approximately equals $t_s=15.3$.}\label{fig:exp}
\end{center}
\end{figure}
See a particular situation from this class represented in Figure \ref{fig:exp}.

\end{ex}

\begin{ex}\label{ex:decay}
Decaying exponential potentials can be considered as well. For
instance for $V(\phi)=(1-e^{-\alpha\sqrt{\phi^2+\gamma^2}})^2$ the
function $u(\phi)$ behaves like $e^{-\alpha|\phi|}$, so it goes to
zero at $\pm\infty$.
\begin{figure}
\begin{center}
\psfull \epsfig{file=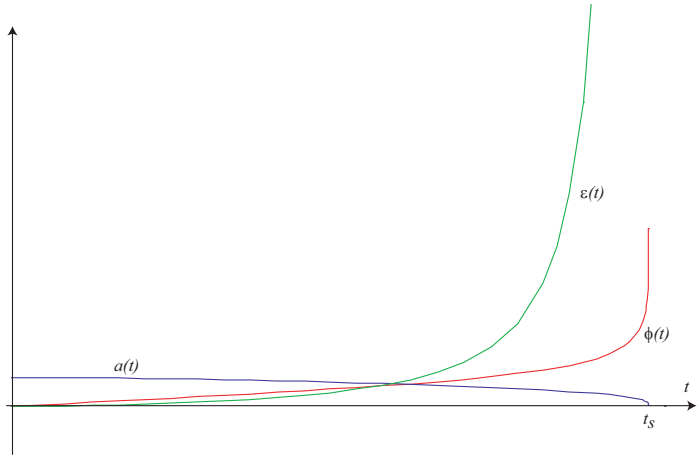, height=6cm}
\caption{Same as Figure \ref{fig:poly}, where now $V(\phi)=2(1-e^{-\sqrt{x^2+1}})^2-\tfrac16 e^{-(x+4)^2}-\tfrac12 e^{-(x-4)^2}+1$.
The initial conditions are $\phi_0=0,\,\dot\phi_0=\sqrt{-2V(0)},\,a_0=1$. The time of collapse $t_s$ approximately equals $t_s=1.8$.}\label{fig:decay}
\end{center}
\end{figure}
Figure \ref{fig:decay} represent a situation where some
corrections terms have been added in order to obtain a potential
with more than one critical points. Even choosing initial data
such that $\epsilon_0=0$, the solution diverges in a finite amount
of time, as showed in general in Subsection \ref{sec:sottosoglia}.
\end{ex}

\section{Gravitational collapse models}

In what follows, we construct models of collapsing objects
composed by homogeneous scalar fields. To achieve this goal we
must match a collapsing solution (considered as the interior
solution in matter) with an exterior spacetime. The natural choice
for the exterior is the so--called \emph{generalized Vaidya
solution}
\begin{equation}\label{eq:Va}
\text ds_{\mathrm{ext}}^2=-\left(1-\frac{2M(U,Y)}Y\right)\,\mathrm dU^2-2\,\mathrm
dY\,\mathrm dU + Y^2\,\mathrm d\Omega^2,
\end{equation}
where $M$ is an arbitrary (positive) function (we refer to
\cite{ww} for a detailed physical discussion of this spacetime,
which is essentially the spacetime generated by a radiating
fluid). The matching is performed along a hypersurface $\Sigma$
which, in terms of a spherical system of coordinates for the
interior, has the simple form  $r=r_b=const$. The Israel junction
conditions at the matching hypersurface read as follows (see
\cite[Proposition 4.1 and Remark 4.2]{hsf}):
\begin{align}
&M(U(t),Y(t))=\frac12 r_b^2 a(t)\dot a(t)^2,\label{eq:mass-match1}\\
&\frac{\partial M}{\partial Y}(U(t),Y(t))=\frac12 r_b^3(\dot a(t)^2+2 a(t)\ddot a(t)),\label{eq:mass-match2}
\end{align}
where the functions $(Y(t),U(t))$ satisfy
\begin{equation}\label{eq:cond}
Y(t)=r_b a(t),\qquad\frac{\text d U}{\text d
t}(t)=\frac1{1+\dot a(t) r_b}.
\end{equation}
The two equations \eqref{eq:mass-match1}--\eqref{eq:mass-match2}
are equivalent to require that Misner--Sharp mass $M$ is continuous and $\tfrac{\partial M}{\partial U}=0$
on the junction hypersurface $\Sigma$.

The endstate of the
collapse of these ``homogeneous scalar field stars'' is analyzed in
the following theorem.

\begin{teo}\label{thm:endstate}
Except at most for a measure zero set of initial data,
the scalar field model matched with the generalized Vaidya solution
\eqref{eq:Va} collapses to a black hole.
\end{teo}
\begin{proof}
We will make use of the result in \cite[Theorem 5.2]{hsf}, that
states that if, in the approach to the singularity, the quantity
$\dot a^2=a^2\epsilon$ is bounded, then the apparent horizon
cannot form, and the singularity is naked. It is indeed easy to
check that the equation of the apparent horizon for the metric
\eqref{eq:g} is given by $r^2\dot a(t)^2=1$. If $\dot a$ is bounded,
one can choose the junction surface $r=r_b$ sufficiently small
such that $\dot a^2(t)<\tfrac1{r^2}$, $\forall
(t,r)\in[0,t_s]\times[0,r_b]$, and so $\left(1-\tfrac{2M}R\right)$
is bounded away from zero near the singularity. As a consequence,
one can find in the exterior portion of the spacetime
\eqref{eq:Va}, null radial geodesics which meet the singularity in
the past, and therefore the singularity is naked. Otherwise, if
$\dot a^2$ is unbounded, the trapped region forms and the collapse
ends into a black hole. Now, we have proved in Lemma
\ref{thm:perla} the existence of
$\rho_\infty=\lim_{\phi\to\infty}\rho(\phi)$ (which is also equal
to $\lim_{t\to t_s^-}\rho(t)$), and that, actually, $\rho_\infty$
vanishes except for a zero--measured set of initial data. On the
other end, using \eqref{eq:original-efe2}, we get
\begin{equation}\label{eq:efe-nu2}
\ddot a=-\frac{\dot a^2}a\left(\frac{2-\rho}{1+\rho}\right),
\end{equation}
and therefore $\ddot a\le-\tfrac{\dot a^2}a$ eventually holds in a
left neighborhood of $t_s$, say $[t_0,t_s[$, where $\dot a$ is
decreasing. It follows that $\forall t\in]t_0,t_s[$, $\dot a(t)\le
\dot a(t_0)$, which is negative by hypotheses. Then
\[
\dot a(t)-\dot a(t_0)
=\int_{t_0}^t\ddot a(\tau)\,\mathrm d\tau\le-\int_{t_0}^t\frac{\dot a(\tau)^2}{a(\tau)}\,\mathrm d\tau\le -\dot a(t_0)\int_{t_0}^t\frac{\dot a(\tau)}{a(\tau)}\,\mathrm d\tau,
\]
that diverges to $-\infty$ by \eqref{eq:a-van}. Then $\dot a(t)$ is unbounded.
\end{proof}

\subsection{Examples of non-generic situations}\label{rem:last}
In the crucial Lemma \ref{thm:perla} the function $\rho$ has been shown to vanish, in the late time behavior, for almost every choice of the initial data.
Now we want to reconsider the non generic situation where
$\lim_{t\to{t_s^-}}{y(\phi(t))}{u(\phi(t))}=1$.
Actually, if the potential $V(\phi)$ satisfies the condition
$$\liminf_{\phi\to\pm\infty}|u(\phi)|> 0,$$
the argument of Proposition \ref{thm:sing} may be used also in this (non generic case),
to prove the formation of the singularity in a finite amount of time.
Indeed, in this case $y$ is bounded above, and then $\rho$ is, so that one can use \eqref{eq:dotenergy} to find a constant $\kappa>0$ such that $\dot\epsilon>\kappa\epsilon^{3/2}$ in a left neighborhood of $\sup\I$.


Moreover, recalling that
$\rho=y^2-1$, and equation \eqref{eq:efe-nu2}, if
$$\liminf_{\phi\to\pm\infty}|u(\phi)|>\tfrac{\sqrt 3}3,$$
the
argument in the proof of Theorem \ref{thm:endstate} also applies
here, and the collapse ends into a black hole. In particular,
observe that in this case $2-\rho$ is bounded away from zero, and
positive. On the other side, if
$$\limsup_{\phi\to\pm\infty}|u(\phi)|<\tfrac{\sqrt 3}3,$$
then
$\ddot a$ is eventually positive, and this implies that $\dot a$
is bounded, so that the formation of the apparent horizon is
forbidden. As a consequence, the collapse, in this non--generic situation, ends into a naked
singularity.

\begin{ex}\label{ex:exp2}
The potential in exponential form of Example \ref{ex:exp} is such that $u(\phi)$ behaves like a positive constant $\lambda$, and then it always gives rise to a black hole.
The -- non generical -- data producing the solution such that $\lim_{\phi\to\infty}\rho(\phi)=\tfrac1{\lambda^2}-1$
are discussed in \cite{hsf, joshi};  a choice of $\lambda<\tfrac{\sqrt 3}3$ makes the limit $\rho_\infty$ greater than 2, and the apparent horizon cannot form, resulting in a naked singularity.
\end{ex}

\begin{ex}\label{ex:poly2}
Let $V$ the potential given by the polynomial in Example
\ref{ex:poly}. We already know that almost every choice of initial
data forms a singularity, which happens to be a black hole. Note
that, in the general case of a polynomial potential with leading
term $\lambda^2 \phi^{2n}$, with $n>2$, it can be shown that
\emph{every} initial data gives rise to a singularity, although
the function $u(\phi)$ diverges and therefore the above arguments cannot be applied.
Indeed, if $y(\phi)$ goes like $u(\phi)^{-1}$ -- that in this case
behaves like $\tfrac n{\sqrt 3\phi}$, then using \eqref{eq:rho} it
is easily seen that $\dot\phi\cong \phi^{n-1}$, since
$$
\frac{3\phi^2}{n^2}\cong y^2=1+\rho\cong 1+\frac{2\lambda^2\phi^{2n}}{\dot\phi^2},
$$
and then if $n>2$ the
solution must diverge in a finite amount of
time. The non generical data yielding this particular situation, moreover, forbid apparent horizon formation, and so the resulting singularity is naked.
\end{ex}

\section{Discussion and Conclusions}

We have discussed here the qualitative behavior of the solutions
of the Einstein field equations with homogeneous scalar fields
sources in dependence of the choice of the self-interacting
potential. In the case of cosmological models, our results widely
extend the recent relevant results obtained by Rendall
(\cite{r1,r2}) and by Miritzis (\cite{mir1}). Using asymptotic
analysis, Rendall has been able to show that scalar fields exhibit
an oscillating behavior if the potential is quadratic (that is, if
the field is free but not massless) or if $V$ is of the form
$a\phi^2 +O(\phi^3)$, while Miritzis has classified the limiting
behaviour of solutions for a large class of non-negative scalar
field potentials. In both these cases however, global existence is
guaranteed. Our approach here extends the "spectrum" of available
potential, and in doing so includes the cases in which
there is no global in time evolution. This allows us to treat in a
unified manner also the case of gravitational collapse, in which a
singularity is always formed in the future. Matching these
solutions with a Vaidya "radiating star" exterior we obtained
models of gravitational collapse which can be viewed as the
scalar-field generalization of the Oppenheimer-Snyder collapse
model, in which a dust homogeneous universe is matched with a
Schwarzschild solution (the Schwarzschild solution can actually be
seen as a special case of Vaidya).

The Oppenheimer-Snyder model,as is well known, describes the
formation of a covered singularity, i.e. a blackhole (it is
actually the first model of blackhole formation ever discovered).
The same occurs here: indeed we show that homogeneous scalr field
collapse generically forms a blackhole. The examples of naked
singularities which were found in recent papers \cite{hsf,joshi},
turn out to correspond to very special cases which,
mathematically, are not generic. Therefore, our results here show
that naked singularities are not generic in homogeneous,
self-interacting scalar field collapse, at least for the
considered (but wide) class of physically relevant potentials.
Non-genericity is already well known for \emph{non
self-interacting} (i.e. $V(\phi)=0$) spherically symmetric scalar
fields. Whether this result can actually be shown to hold also in
the much more difficult case of both inhomogeneous and
self-interacting scalar fields remains an open problem.

\appendix

\section{An existence/uniqueness
theorem for a kind of singular ODE}\label{sec:proof}

In this section we prove a result of existence, uniqueness, and instability of
solution for a particular kind of ordinary differential equation
of first order to which the standard theory does not apply
straightforwardly. To the best of our knowledge this result,
needed in the proof of Lemma \ref{thm:perla}, although relatively
simple, does not appear in the literature.

\begin{teo}\label{thm:ode}
Let us consider the Cauchy problem
$$
h(s)\dot z(s)=f(s,z(s))+g(s),\qquad z(0)=0.
$$
and $\beta$ is a positive constant such that the following conditions hold:
\begin{enumerate}
\item $h\in{\mathcal C}^0([0,\beta],\R)$ such that $h(0)=0$, $h(s)>0$ in $]0,\beta]$, and  $h(s)^{-1}$ is not integrable in $]0,\beta]$;
\item $g\in{\mathcal C}^0([0,\beta],\R)$ such that $g(0)=0$;
\item $f\in{\mathcal C}^{0,1}([0,\beta]\times\R,\R)$ such that $f(s,0)\equiv 0$, $\tfrac{\partial f}{\partial z}(0,0)< 0$, and
\item $\exists\rho>0$ such that $\tfrac{\partial f}{\partial z}(s,z)$ is  uniformly Lipschitz continuous with respect to $z$ in $[0,\beta]\times[-\rho,\rho]$, that is $\exists L>0$ such that, if $\vert z_1\vert, \vert z_2\vert\le\rho$, $s\in [0,\beta]$, then $\vert \tfrac{\partial f}{\partial z}(s,z_1)-\tfrac{\partial f}{\partial z}(s,z_2)\vert\le L \vert z_1-z_2\vert$.
\end{enumerate}

Then, there exists $\alpha<\beta$ such that the above Cauchy problem admits a unique solution $z(s)$ in $[0,\alpha]$. Moreover, this solution is the only function solving the differential equation $h(s)\dot z(s)=f(s,z(s))+g(s)$ with the further property $\liminf_{s\to 0}\vert z(s)\vert=0$.

\end{teo}

\begin{proof}
Let us set $\alpha<\beta$ free to be determined later, and let ${\mathcal X}_\alpha$ be the space
$$
{\mathcal X}_\alpha=\{z\in{\mathcal C}^0([0,\alpha])\cap{\mathcal C}^1(]0,\alpha])\,:\,z(0)=0,\,\lim_{s\to 0^+}h(s)\dot z(s)=0\}.
$$
It can be proved that ${\mathcal X}_\alpha$ is a Banach space, endowed with the norm
$\Vert z\Vert_\alpha=\Vert z\Vert_\infty + \Vert h\dot z\Vert_\infty$.

Let also be
$$
{\mathcal Y}_\alpha=\{\lambda\in{\mathcal C}^0([0,\alpha])\,:\,\lambda(0)=0\}
$$
a (Banach) space endowed with the $L^\infty$--norm and let us consider the functional
$$
{\mathcal F}:{\mathcal X}_\alpha\to{\mathcal Y}_\alpha, \qquad {\mathcal F}(z)(s)=h(z)\dot z(s)-f(s,z(s)).
$$
It is easily verified that $\mathcal F$ is a $\mathcal C^1$ functional, with tangent map at a generic element $z\in\mathcal X_\alpha$ given by
\[
\left(\text d\mathcal F(z)[\xi]\right)(s)= h(s)\dot\xi(s)-\frac{\partial f}{\partial z}(s,z(s))\xi(s),
\]
where $\xi\in  X_\alpha$.
Observing that $g\in{\mathcal Y}_\alpha$,
we want to find $\alpha$ such that the equation
\begin{equation}\label{eq:inv}
\mathcal F(z)=g
\end{equation}
has a unique solution $z \in  X_\alpha$.
To this aim, we will exploit an Inverse Function scheme, and we will prove that $\mathcal F$ is a local homeomorphism from
a neighborhood of $z_0\equiv 0$ in $\mathcal X_\alpha$ onto a neighborhood of $\mathcal F(z_0)\equiv 0$ in $\mathcal Y_\alpha$, that includes $g$.
This will be done using neighborhoods with radius independent of $\alpha$ and this will be crucial to obtain the uniqueness result.

In the following we review the classic scheme (see for instance \cite{Berg}) for reader's convenience. Let $R(z)=\mathcal F(z)-\mathcal F(0)-\text d\mathcal F(0)[z]$; then \eqref{eq:inv} is equivalent to find $z$ such that $R(z)+\text d\mathcal F(0)[z]=g$, and therefore, if $d\mathcal F(0)$ is invertible, to prove the existence of a unique fixed point of the application $T$ on $\mathcal X_\alpha$,
\begin{equation}\label{eq:T}
T(z)=(d\mathcal F(0))^{-1}[g-R(0,z)].
\end{equation}
We first show that $T$ is a contraction map from the ball $B(0,\delta)\subseteq\mathcal X_\alpha$ in itself, provided that $\delta$
and $\Vert g\Vert_\infty$ are
sufficiently small (independently by $\alpha$). The following facts must be proven to this aim:
\begin{enumerate}
\item there exists a constant $M$, independent on $\alpha$,
such that, $\forall z_1,z_2\in\mathcal X_\alpha$ with $\Vert z_1\Vert_\alpha,\Vert z_2\Vert_\alpha\le 1$, it is
\begin{equation}\label{eq:F-lip}
\Vert\text d\mathcal F(z_1)-\text d\mathcal F(z_2)\Vert
\le M \Vert z_1-z_2\Vert_\alpha
\end{equation}
(the norm on the left hand side refers to the space of linear applications from $\mathcal X_\alpha$ to $\mathcal Y_\alpha$).
\item there exists a constant $C$, independent on $\alpha$, such that
\begin{equation}\label{eq:Fprime}
\Vert \text d\mathcal F(0)^{-1}\Vert \le C
\end{equation}
(here the norm refers to the space of linear applications from $\mathcal Y_\alpha$ to $\mathcal X_\alpha$).
\end{enumerate}
If the above facts hold, given $z_1,z_2\in \mathcal X_\alpha$, then
\begin{multline*}
\text d\mathcal F(0)[T(z_1)-T(z_2)]\\ =R(z_2)-R(z_1)=\mathcal F(z_2)-\mathcal F(z_1)-\text d\mathcal F(0)[z_2-z_1]\\ =\int_0^1 \left(\text d\mathcal F(t z_2+(1-t)z_1)-\text d\mathcal F(0)\right)[z_2-z_1]\,\text dt,
\end{multline*}
hence, if in addition $z_1,z_2\in B(0,\delta)$,
\begin{multline*}
\Vert T(z_1)-T(z_2)\Vert_\alpha\le \\ \Vert (\text d\mathcal F(0))^{-1}\Vert \left(\int_0^1
\Vert \left(\text d\mathcal F(t z_2+(1-t)z_1)-\text d\mathcal F(0)\right)\Vert \,\text dt\right)\,\Vert z_1-z_2\Vert_\infty\le 2 M\,C\,\delta
\Vert z_1-z_2\Vert_\alpha.
\end{multline*}
and therefore:
\begin{itemize}
\item $T$ is a contraction, taking $\delta$ such that $K:=2MC\delta<1$;
\item since $\Vert T(z)\Vert_\alpha\le \Vert T(z)-T(0)\Vert_\alpha+\Vert T(0)\Vert_\alpha\le K \Vert z\Vert_\alpha+ \Vert (d\mathcal F(0))^{-1}[g]
\Vert_\alpha\le K \Vert z\Vert_\alpha+C\Vert g\Vert_\infty $, then $T$ maps $B(0,\delta)$ in itself,
provided that $\Vert g\Vert_\infty\le C^{-1}(1-K)\delta$.
\end{itemize}
Observe that the first of these two facts fixes the value of $\delta$, whilst the last inequality holds choosing $\alpha$ -- free so far -- small enough.
This is one of the reasons why the constants $M$ and $C$ must be independent on $\alpha$.
Then $T$ admits a unique fixed point on $B(0,\delta)$, which is the solution to our problem \eqref{eq:inv} on the interval $[0,\alpha]$. In other words, the function $\mathcal F$ is a local homeomorphism from $B(0,\delta)\subseteq\mathcal X_\alpha$ to $B(0,C^{-1}(1-K)\delta)\subseteq\mathcal Y_\alpha$.

To see that the solution is \emph{globally} unique on $\mathcal X_\alpha$, let us argue
as follows. Suppose $\bar z\in\mathcal X_\alpha\setminus B(0,\delta)$ $(\bar z\ne z)$
solves the problem, and let $\alpha_1\le\bar\alpha$ sufficiently small such that $\Vert
\bar z\Vert_{\alpha_1}\le\delta$. Then, observing $\Vert
g\Vert_{L^\infty([0,\alpha_1])}\le\Vert g\Vert_{L^\infty([0,\alpha])}\le
C^{-1}(1-K)\delta$, and recalling that estimates \eqref{eq:F-lip}--\eqref{eq:Fprime} do not depend on $\alpha$, one can argue as before to find that $\bar z\vert_{[0,\alpha_1]}$ is the unique element of $B(0,\delta)\subseteq\mathcal X_{\alpha_1}$ mapped into
$g\vert_{[0,\alpha_1]}\in B(0,C^{-1}(1-K)\delta)\subseteq\mathcal Y_{\alpha_1}$. But, of course, $\Vert z\Vert_{\alpha_1}\le\Vert z\Vert_\alpha<\delta$, so $z$ and $\bar z$ coincide on $[0,\alpha_1]$, and therefore on all $[0,\alpha]$.

Therefore, to complete the proof of the existence and uniqueness for the given Cauchy problem, just \eqref{eq:F-lip}--\eqref{eq:Fprime} are to be proven. The first equation is a consequence of local
uniform Lipschitz continuity of $f_{,z}$. The second one needs some more care: taken $\lambda\in\mathcal Y_\alpha$, we must consider the Cauchy problem
\begin{equation}\label{eq:linCauchy}
h(s)\dot\xi(s)= \ell(s) \xi(s)+\lambda(s),\qquad\xi(0)=0,
\end{equation}
where $\ell(s):=\tfrac{\partial f}{\partial z}(s,0)$, that without loss of generality we can suppose negative, $\forall s\in [0,\beta]$. First, it is easily seen that \eqref{eq:linCauchy} admits the unique solution $\xi\in\mathcal X_\alpha$,
\[
\xi(s)=e^{-\int_s^\alpha{\ell(t)}{h(t)^{-1}}\,\text dt}\int_0^s\frac{\lambda(t)}{h(t)}
e^{\int_t^\alpha{\ell(\tau)}{h(\tau)^{-1}}\,\text d\tau}\,\text dt.
\]
Then $\Vert(\text d\mathcal F(0))^{-1}[\lambda]\Vert_\alpha=\Vert\xi\Vert_\alpha\le(1+\ell_1)\Vert\xi\Vert_\infty+\Vert\lambda\Vert_\infty$,
where $\ell_1=\Vert\ell\Vert_{L^\infty{[0,\beta]}}$. Moreover, called $\ell_0=\sup_{[0,\beta]}\ell(s)<0$, it is easily seen that
$\Vert\xi\Vert_\infty\le-\tfrac1{\ell_0}\Vert\lambda\Vert_\infty$, so it suffices to choose $C=1-\tfrac{\ell_1+1}{\ell_0}$, and \eqref{eq:Fprime} is proven.

To prove last claim of the Theorem, let us suppose that $w(s)$ is a function defined in $[0,\alpha]$ such that $h(s)\dot w(s)=f(s,w(s))+g(s)$, and that $s_k$ is an infinitesimal and monotonically decreasing sequence such that $w(s_k)\to 0$ as $k\to\infty$. We want to prove that $w=z$, and therefore, it will suffice to show that $\lim_{s\to 0} w(s)=0$.

First of all, observe that from the hypotheses, the equation $f(s,z)+g(s)=0$ defines a continuous function $\zeta(s):[0,\delta]\to\R$, such that $\zeta(0)=0$.
In particular, since $\frac{\partial f}{\partial z}(0,0)<0$, $\exists\rho>0$ such that, in the rectangle $[0,\delta]\times[-\rho,\rho]$, it must be
$f(s,z)+g(s)<0$ (resp.: $>0$) if $z>\zeta(s)$ (resp.:$z<\zeta(s)$).

Let us now argue by contradiction, supposing the existence of an infinitesimal sequence $\sigma_k$, that can be chosen with the property $\sigma_k<s_k$, such that $\vert w(\sigma_k)\vert>\theta$ for some given constant $\theta$. Now, up to taking a smaller constant $\delta$, then $\vert\zeta(s)\vert<\tfrac\theta 2,\,\forall s\in[0,\delta]$.
Then, for $k$ sufficiently large, $\vert w(\sigma_k)\vert\ge\theta>\tfrac\theta 2\ge\sup_{[0,\delta]}\vert\zeta(s)\vert$. Therefore, $\forall s<\sigma_k$, $\vert w(s)\vert>\vert w(\sigma_k)\vert$, which is a contradiction since $w(s_k)\to 0$. Then $\lim_{s\to 0}w(s)=0$, and the proof is complete.

\end{proof}

\section{Local existence/uniqueness of solutions with initial zero--energy}\label{sec:eps0}

\begin{lem}\label{lem:eps0}
Let $\phi_0,\,v_0$ such that $v_0^2+2V(\phi_0)=0$. Then, there exists $t_*>0$ such that the Cauchy problem
\begin{equation}\label{eq:eps0}
\begin{cases}
&\ddot\phi(t)=-V'(\phi(t))+\sqrt{3(\dot\phi(t)^2+2V(\phi(t)))}\,\dot\phi(t),\\
&\phi(0)=\phi_0,\\
&\dot\phi_0=v_0,
\end{cases}
\end{equation}
has a unique solution $\phi(t)$ defined in $[0,t_*]$ with the property
\begin{equation}\label{selezione}
\epsilon(t)=3\left(\int_0^t\dot\phi(s)^2\,\text ds\right)^2,\,
\forall t\in]0,t_*].
\end{equation}

Moreover if $(\phi_{0,m},v_{0,m}) \rightarrow (\phi_0,v_0)$, $(v_{0,m})^2 + 2V(\phi_{0,m}) = 0$ and $\phi_m$ is the solution of \eqref{eq:eps0}
with initial data
$(\phi_{0,m},v_{0,m})$ satisfying condition \eqref{selezione}, it is $\phi_m \rightarrow \phi$ with respect to the $C^2$-norm in the interval $[0,t_*]$.
\end{lem}
\begin{proof}
Let us consider the "penalized" problem
\begin{equation}\label{eq:pen}
\begin{cases}
&\ddot\phi(t)=-V'(\phi(t))+\sqrt{3(\dot\phi(t)^2+2V(\phi(t))+\frac1{n^2})}\,\dot\phi(t),\\
&\phi(0)=\phi_0,\\
&\dot\phi_0=v_0,
\end{cases}
\end{equation}
that has a unique local solution $\phi_n$. If $\phi_n$ is not defined $\forall t\ge 0$, let $\I_n$ be the set
\[
\I_n=\{t\in\R\,:\,|\phi_n(s)| \leq |\phi_0| +1,\,|\dot\phi_n(s)|\le|v_0|+1,\,\forall s\ge t\}.
\]
Of course, $\I_n\ne\emptyset$ and, called $t_n=\sup\I_n$, if $t_n$ is finite, then $|\dot\phi_n(t_n)|=|v_0|+1$,
or
$|\phi_n(t_n)| = |\phi_0| +1$. Now assume $|\phi_n(t_n)| = |\phi_0| +1$. Then
\[
1 = \vert \phi_n(t_n)-\phi_0\vert\le \int_0^{t_n}\vert\dot\phi_n(s)\vert\,\text ds\le(|v_0|+1)t_n.
\]
Analogously if  $|\dot\phi_n(t_n)|=|v_0|+1$ we have
\[
1 = \vert \dot \phi_n(t_n)- v_0\vert\le \int_0^{t_n}\vert\ddot\phi_n(s)\vert\,\text ds.
\]
Since $|\phi_n(t)| \leq |\phi_0| +1$ and $|\dot\phi_n(s)|\le|v_0|+1$ for all $t \in [0,t_n]$, and $\phi_n$ solves \eqref{eq:pen}, we see that there exists
$C$ independent of $n$ such that $\vert\ddot \phi(t)\vert \leq C$ for all $t \in [0,t_n]$. Therefore in this second case we obtain $1 \leq Ct_n$.

Then $t_*:=\inf_n t_n>0$ (we set $t_*=1$ if $t_n=+\infty$ $\forall n$). Moreover $\vert\ddot \phi_n\vert$ is uniformly bounded in $[0,t_*]$, then up to subsequences, there exists a $\mathcal C^1$ function $\phi(t)$, solution of \eqref{eq:eps0}, such that $\phi_n\to\phi$ and $\dot\phi_n\to\dot\phi$ uniformly on $[0,t_*]$.

Now, consider $\epsilon_n(t):=\dot\phi_n(t)^2+2V(\phi_n(t))+\tfrac1{n^2}$. We have
\begin{equation}\label{energiaappr}
\dot\epsilon_n(t)=2\sqrt 3\sqrt{\epsilon_n(t)}\dot\phi_n(t)^2.
\end{equation}
Then $\epsilon$ is not decreasing, while $\epsilon(0) = \frac1n$. Then is uniformly bounded away from zero and therefore by
\eqref{energiaappr}, dividing by $\sqrt{\epsilon_n}$ and integrating gives
$\sqrt{\epsilon_n(t)}=\tfrac1{n^2}+\sqrt 3\int_0^t \dot\phi_n(s)^2\,\text ds$.
Therefore passing to the limit in $n$ we obtain $\epsilon(t)=\dot\phi(t)^2+2V(\phi(t))=3\left(\int_0^t \dot\phi(s)^2\,\text ds\right)^2$ for all
$t \in [0,t_*]$ obtaining the proof of the existence of a solution.

The uniqueness of such a solution can be obtained by a contradiction argument. Assuming $\phi$ and $\psi$ solutions, and called $\theta=\phi-\psi$, one can obtain, using \eqref{eq:eps0}, the estimate
\[
\vert\dot\theta(t)\vert\le K_1\int_0^t\vert\theta(s)\vert\,\text ds+ K_2\int_0^t\vert\dot\theta(s)\vert\,\text ds,
\]
for suitable constants $K_1,K_2$. Setting $\rho(t)=\vert\theta(t)\vert+\vert\dot\theta(t)\vert$, and observing that
$\rho(0)=0$, it is not hard to get the estimate
$\rho(t)\le (K_1+K_2+1)\int_0^t\rho(s)\,\text ds$, and then $\rho\equiv 0$ from Gronwall's inequality.

Finally using Gronwall's Lemma as above we obtain also the continuity with respect to the initial data.
\end{proof}

\begin{rem}\label{rem:eps0rev}
Reversing time direction in the above discussed problem \eqref{eq:eps0} yields a results of genericity for expanding solutions such that the energy $\epsilon(t)$ vanishes at some finite time $T$.
\end{rem}


\begin{thebibliography}{00}



\bibitem{Berg} M.S. Berger,  \emph{Nonlinearity and Functional Analysis}, Academic Press: New York, 1977.

\bibitem{c1} D. Christodoulou,
Ann. Math. {\bf 140} 607 (1994).

\bibitem{c2} D. Christodoulou,
Ann. Math. {\bf 149} 183 (1999).



\bibitem{ns} R. Giamb\`o, F. Giannoni, G. Magli, P. Piccione, Comm. Math. Phys. \textbf{235(3)}  545-563 (2003)

\bibitem{hsf} R. Giamb\`o, Class. Quantum Grav. \textbf{22} (2005) 1-11


\bibitem{jmp} R. Giamb\`o, F. Giannoni, G. Magli, J. Math. Phys., \textbf{47} 112505 (2006)

\bibitem{joshi} R. Goswami, P. S. Joshi, {\tt gr-qc/0410144}



\bibitem{Global}  P. S. Joshi,
{\it Global aspects in gravitation and cosmology}, (Clarendon
press, Oxford, 1993).

\bibitem{j-rev} P. S. Joshi, Modern Phys. Lett. A \textbf{17}  1067--1079 (2002)

\bibitem{mir1} J. Miritzis,
Class. Quantum Grav.  \textbf{20}  (2003),  no. 14, 2981--2990

\bibitem{mir2} J. Miritzis, J. Math. Phys. \textbf{44} (2003) 3900-3910

\bibitem{mir3} J. Miritzis, J. Math. Phys. \textbf{46} (2005) 082502



\bibitem{perko} L. Perko, \emph{Differential Equations and Dynamical Systems}, Springer--Verlag, New York, 1991

\bibitem{r1} A.D. Rendall Class.Quant.Grav. \textbf{21} (2004) 2445-2454

\bibitem{r2} A.D. Rendall Class.Quant.Grav. \textbf{24} (2007) 667-678
\bibitem{st1} C. Rubano, J. D. Barrow, Phys.Rev. \textbf{D64} (2001) 127301


\bibitem{st2} C. Rubano, P. Scudellaro Gen.Rel.Grav. \textbf{34} (2002) 307-328



\bibitem{ww} A. Wang and Y. Wu, 1999 Gen. Rel. Grav. \textbf{31} 107


\end{thebibliography}
\end{document}